\documentclass[12pt]{article}
\usepackage{amssymb,latexsym}
\def\ekv#1#2{\begeq\label{#1}#2\endeq}
\def\eekv#1#2#3{\begin{eqnarray}\label{#1}#2\\ #3\nonumber\end{eqnarray}}

\def\aby{arbitrary}

\def\an{analytic}
\def\asy{asymptotic}
\def\bdd{bounded}

\def\coef{coefficient}

\def\canform{canonical transformation}

\def\ctf{canonical transformation}

\def\dop{differential operator}
\def\ef{eigenfunction}
\def\ev{eigenvalue}
\def\e{equation}
\def\fu{function}
\def\fy{family}
\def\F{Fourier}
\def\fop{Fourier integral operator}

\def\hol{holomorphic}
\def\indep{independent}

\def\mfld{manifold}
\def\ml{microlocal}
\def\neigh{neighborhood}
\def\nondeg{non-degenerate}
\def\op{operator}

\def\pb{problem}

\def\pert{perturbation}

\def\pol{polynomial}

\def\pop{pseudodifferential operator}

\def\res{resonance}

\def\sa{selfadjoint}

\def\sop{Schr{\"o}dinger operator}
\def\st{strictly}

\def\sufly{sufficiently}

\def\Th{Theorem}
\def\th{theorem}

\def\uf{uniform}

\def\vf{vector field}
\def\wrt{with respect to}
\def\Re{{\rm Re\,}}
\def\Im{{\rm Im\,}}
\title{Perturbations of \sa{} \op{}s with periodic classical 
flow}
\author{Johannes Sj{\"o}strand\\Centre de Math{\'e}matiques, Ecole Polytechnique\\FR 
91120 Palaiseau, France\footnote{UMR 7640, CNRS}\\johannes@math.polytechnique.fr}
\date{}
\def\wrtext#1{\relax\ifmmode{\leavevmode\hbox{#1}}\else{#1}\fi}

\def\begeq{\begin{equation}}
\def\endeq{\end{equation}}

\newcommand{\eps}{\epsilon}
\def\part#1{\frac{\partial}{\partial #1}}

\newcommand{\real}{\mbox{\bf R}}

\newcommand{\z}{\mbox{\bf Z}}

\renewcommand{\Re}{\mbox{\rm Re\,}}
\renewcommand{\Im}{\mbox{\rm Im\,}}

\renewcommand{\exp}{\mbox{\rm exp\,}}

\newtheorem{dref}{Definition}[section]

\newtheorem{theo}[dref]{Theorem}
\newtheorem{prop}[dref]{Proposition}

\begin{document}

\maketitle
\begin{abstract}
We consider non-\sa{} perturbations 
of a self-adjoint $h$-pseudo\-differ\-ential \op{} in dimension 2. 
In the present work we treat 
the case when the classical flow of the unperturbed part is periodic and 
the strength $\epsilon $ of the perturbation satisfies $h^{\delta_0} <\epsilon 
\le \epsilon _0$ for some $\delta_0\in ]0,1/2[$ and a \sufly{} small 
$\epsilon _0>0$. We get a 
complete \asy{} description of all \ev{}s in certain rectangles 
$[-1/C,1/C]+i\epsilon [F_0-1/C,F_0+1/C]$. In particular we are able to 
treat the case when $\epsilon >0$ is small but \indep{} of $h$. 
\end{abstract}

\section{Introduction}\label{section0}
\setcounter{equation}{0}

\par This paper is a continuation of  \cite{MeSj}, where A. Melin and the author 
observed that for a wide
and stable class of non-\sa{} \op{}s in dimension 2 and in the
semi-classical limit ($h\to 0$), it is possible to describe all \ev{}s
individually in an $h$-\indep{} domain in ${\bf C}$, by means of a
Bohr-Sommerfeld quantization condition. Notice that the corresponding conclusion 
in the \sa{} case seems to be
possible only in dimension 1, or in higher dimensions under strong (and unstable) assumptions
of complete integrability. In \cite{MeSj} we exploited the absence of small 
denominators to get a geometric analogue of the KAM-\th{} via methods of 
non-linear Cauchy--Riemann \e{}s and got a corresponding result at the 
level of \op{}s 

\par In the present work we make another step by studying small 
perturbations, roughly of the form $P+i\epsilon Q$, of a \sa{}  
$h$-\pop{} $P$ whose associated classical flow is periodic. We will here 
be particularly interested in the case of a small but fixed $\epsilon $, 
but our methods allow us to let $\epsilon $ vary in an interval 
$[h^{\delta_0} ,\epsilon _0]$ where $\epsilon _0>0$ is \sufly{} small 
and $\delta _0\in ]0,1/2[$ is \aby{}. 

\par From the point of view of applications, it is clear that even smaller 
perturbations are of a considerable interest and as another step, Hitrik and 
the author \cite{HiSj} studied the same problem 
as in the present paper, but in the parameter range $h\ll \epsilon 
\le h^\delta $ for every fixed $\delta >0$. When the subprincipal 
symbol vanishes we could even treat the range $h^2\ll \epsilon \le 
h^\delta $. Actually with M. Hitrik, we plan a whole 
series of works devoted to small 
perturbations of non-\sa{} \op{}s in two dimensions. Among other things we plan to treat the 
case when the classical flow of the unperturbed \op{} admits certain invariant 
Lagrangian torii with a diophantine condition. (Another work (\cite{Sj2}) 
deals with \res{}s generated 
by a closed hyperbolic trajectory and can be viewed as descendant of the 
pioneering work of M. Ikawa \cite{Ik} about scattering poles for 
two strictly convex obstacles.) 

\par The methods in \cite{HiSj} are partly 
more traditional and rely on reduction by averaging to a one-dimensional 
problem in the spirit of \cite{We, Co, BoGu, HeRo, Iv}. Such a 
reduction does not seem possible here and the problem 
remains 2-dimensional. In general, we have been motivated by recent progress
around the damped wave-equation (\cite{Le}, \cite{AsLe}, \cite{Sj3}, 
\cite{Hi}),  as well as
the problem of barrier top \res{}s for the semi-classical \sop{}
(\cite{KaKe})  where more complete results than the corresponding ones for
\ev{}s of potential wells (\cite{Sj4}, \cite{BaGrPa}, \cite{Po}) seem possible. 
Eventually we also hope to apply our results (though not specifically the 
ones of the present work) to the distribution of \res{}s for a
\st{} convex obstacle in ${\bf R}^3$. See \cite{SjZw2} and references given
there. In the case of analytic obstacles, much more can probably
be said, especially in dimension 3 (and 2).

\par The present work was undertaken before the start of  \cite{HiSj}, 
but the latter work is now completed, so we can take advantage of many of 
the arguments there, even though the main step here will be quite different. 

\par Let $M$ denote ${\bf R}^2$ or a compact real-analytic \mfld{} of
dimension 2.

\par When $M={\bf R}^2$, let 
\begeq
\label{0.1}
P_\epsilon =P(x,hD_x,\epsilon ;h)
\endeq
be the $h$-Weyl quantization on ${\bf R}^2$ of a symbol $P(x,\xi ,\epsilon
;h)$ depending smoothly on $\epsilon \in{\rm neigh\,}(0,{\bf R})$ with
values in the space of \hol{} \fu{}s of $(x,\xi )$ in a tubular
\neigh{} of ${\bf R}^4$ in ${\bf C}^4$, with 
\begeq\label{0.2}
\vert P(x,\xi ,\epsilon ;h)\vert \le Cm(\Re (x,\xi ))
\endeq 
there. Here $m$ is assumed to be an order \fu{} on ${\bf R}^4$, in the
sense that $m>0$ and 
\begeq\label{0.3}
m(X)\le C_0\langle X-Y\rangle ^{N_0}m(Y),\ X,Y\in{\bf R}^4.
\endeq
We also assume that 
\begeq\label{0.4}
m\ge 1.\endeq
We further assume that 
\begeq\label{0.5}
P(x,\xi ,\epsilon ;h)\sim \sum_{j=0}^\infty  p_{j,\epsilon }(x,\xi)h^j,\ h\to 0, \endeq
in the space of such \fu{}s. We make the ellipticity assumption
\begeq\label{0.6}
\vert p_{0,\epsilon }(x,\xi )\vert \ge {1\over C}m(\Re (x,\xi )),\
\vert (x,\xi )\vert \ge C, \endeq
for some $C>0$.

\par When $M$ is a compact \mfld{}, we let 
\begeq\label{0.7}
P_\epsilon =\sum_{\vert \alpha \vert \le m}a_{\alpha ,\epsilon
}(x;h)(hD_x)^\alpha ,\endeq
be a \dop{} on $M$, such that for every choice of analytic local coordinates,
centered at some point of $M$, $a_{\alpha ,\epsilon }(x;h)$ is a smooth
\fu{} of $\epsilon $ with values in the space of \bdd{} \hol{} \fu{}s
in a complex \neigh{} of $x=0$. We further assume that 
\begeq\label{0.7.5}
a_{\alpha ,\epsilon }(x;h)\sim \sum_{j=0}^\infty  a_{\alpha ,\epsilon
,j}(x)h^j,\ h\to 0, \endeq
in the space of such \fu{}s. The semi-classical principal symbol in
this case is given by
\begeq\label{0.8}
p_{0,\epsilon }(x,\xi )=\sum a_{\alpha ,\epsilon ,0}(x)\xi ^\alpha ,
\endeq
and we make the ellipticity assumption
\begeq\label{0.9}
\vert p_0(x,\xi )\vert \ge {1\over C}\langle \xi \rangle ^m,\ (x,\xi
)\in T^*M,\,\vert \xi \vert \ge C,\endeq
for some large $C>0$. (Here we assume that $M$
has been equipped with some Riemannian metric, so that $\vert \xi \vert
$ and $\langle \xi \rangle =(1+\vert \xi \vert ^2)^{1/2}$ are
well-defined.)

\par Sometimes, we write $p_\epsilon $ for $p_{0,\epsilon }$ and
simply $p$ for $p_{0,0}$. Assume 
\begeq\label{0.10}
P_{\epsilon =0} \hbox{ is formally \sa{}.}
\endeq
In the case when $M$ is compact, we let the under\-lying Hil\-bert
spa\-ce be $L^2(M,\mu (dx))$ for some po\-si\-ti\-ve real-ana\-ly\-tic den\-si\-ty 
$\mu(dx)$ on $M$. 

\par Under these assumptions $P_\epsilon $ will have discrete spectrum
in some fixed \neigh{} of $0\in{\bf C}$, when $h>0,\epsilon \ge 0$
are \sufly{} small, and the spectrum in this region, will be contained
in a band $\vert \Im z\vert \le {\cal O}(\epsilon )$. The purpose of
this work as well as of \cite{HiSj} and later ones in this series, is to give detailed 
\asy{} results about the distribution of individual \ev{}s inside such a band. 

\par Assume for simplicity that (with $p=p_{\epsilon =0}$)
\begeq\label{0.11}
\Gamma _0:=p^{-1}(0)\cap T^*M\hbox{ is connected.}\endeq
Let $H_p=p'_\xi \cdot {\partial \over \partial x}-p'_x\cdot {\partial
\over \partial \xi }$ be the Hamilton field of $p$. In this work, we
will always assume that for $E\in{\rm neigh\,}(0,{\bf R})$:
\begin{eqnarray}\label{0.12}
&&\hbox{The }H_p\hbox{-flow is periodic on }\Gamma _E:=p^{-1}(E)\cap T^*M\hbox{
with}\\ &&\hbox{a period }T(E)>0 \hbox{ depending \an{}ally on }E.\nonumber
\end{eqnarray}
(In Section \ref{SectionGe} we recall how this assumption follows from a 
seemingly weaker one.)
Let $q={1\over i}{({\partial \over \partial \epsilon })}_{\epsilon
=0}p_\epsilon $, so that 
\begeq\label{0.13}
p_\epsilon =p+i\epsilon q+{\cal O}(\epsilon ^2m),\endeq
in the case $M={\bf R}^2$ and $p_\epsilon =p+i\epsilon q+{\cal
O}(\epsilon ^2\langle \xi \rangle ^m)$ in the \mfld{} case. Let 
\begeq\label{0.14}
\langle q\rangle ={1\over T(E)}\int_{-T(E)/2}^{T(E)/2}q\circ \exp (tH_p)\,
dt\hbox{ on }p^{-1}(E)\cap T^*M.\endeq
Notice that $p,\langle q\rangle $
are in involution; $0=H_p\langle q\rangle =:\{ p,\langle q\rangle \} $.
As in \cite{HiSj}, we shall see how to reduce ourselves to the case when
\begeq\label{0.15}
p_\epsilon =p+i\epsilon \langle q\rangle +{\cal O}(\epsilon ^2),\endeq
near $p^{-1}(0)\cap T^*M$. An easy consequence of this is that the
spectrum of $P_\epsilon $ in $\{z\in {\bf C}; \vert \Re z\vert <\delta
\}$ is confined to $]-\delta ,\delta [+i\epsilon ]\langle \Re q\rangle
_{{\rm min},0}-o(1),\langle \Re q\rangle
_{{\rm max},0}+o(1)[$, when $\delta ,\epsilon ,h\to 0$, where $\langle
\Re q\rangle _{{\rm min},0}=\min_{p^{-1}(0)\cap T^*M}\langle \Re
q\rangle $ and similarly for $\langle q\rangle _{{\rm max},0}$. We will
mainly think about the case when $\langle q\rangle $ is real-valued but
we will work under the more general asumption that 
\begeq\label{0.16}
\Im \langle q\rangle \hbox{ is an \an{} \fu{} of }p\hbox{ and }\Re
\langle q\rangle ,\endeq
in the region of $T^*M$, where $\vert p\vert \le 1/{\cal O}(1)$.

\par Let $F_0\in [\langle \Re q\rangle _{{\rm min},0},\langle \Re
q\rangle _{{\rm max},0}]$. The purpose of the present work is to
determine all \ev{}s in a rectangle
\begeq\label{0.17}
]-{1\over {\cal O}(1)},{1\over {\cal O}(1)}[+i\epsilon ]F_0-{1\over
{\cal O}(1)},F_0+{1\over {\cal O}(1)}[, \endeq
for 
\begeq\label{0.18}
h^{\delta_0} \le \epsilon \le \epsilon _0,
\endeq
for $h^{\delta _0}\le \epsilon \le \epsilon _0$ with $0<\delta _0<1/2$ 
and $\epsilon _0$ \sufly{} small but fixed. 
We assume that 
\begeq\label{0.19}
T(0)\hbox{ is the minimal period of every }H_p\hbox{-trajectory in
}\Lambda _{0,F_0},\endeq
where
\begeq\label{0.20}
\Lambda _{0,F_0}:=\{ \rho \in T^*M;\, p(\rho )=0,\, \Re \langle
q\rangle (\rho )=F_0\} ,\endeq
We also assume that 
\begeq\label{0.21}
dp,\, d\Re \langle q\rangle \hbox{ are linearly \indep{} at every
point of }\Lambda _{0,F_0}.\endeq
This implies that every connected component of $\Lambda _{0,F_0}$ is a
two-dimensional Lagrangian torus. For simplicity, we shall assume that
there is only one such component. Notice that in view of (\ref{0.19}), the
space of closed orbits in $p^{-1}(0)\cap T^*M$;
$$\Sigma :=(p^{-1}(0)\cap T^*M)/\sim ,$$
where $\rho \sim \mu $ if $\rho =\exp tH_p \mu $ for some $t\in{\bf
R}$, becomes a 2-dimensional symplectic \mfld{} near the image of
$\Lambda _{0,F_0}$, and (\ref{0.21}) simply means that $\Re \langle q\rangle $, viewed
as a \fu{} on $\Sigma $, has non-vanishing differential.  The image of
$\Lambda _{0,F_0}$ is just a closed curve. 

In \cite{HiSj} (for $\epsilon $ in the range $h\ll \epsilon \le h^\delta $ 
and sometimes $h^2 \ll \epsilon \le h^\delta ,\,\forall \delta >0$) we 
also studied the case 
when $F_0$ is a \nondeg{} extreme valule of $\langle q\rangle $ on 
$\Sigma $. It would be interesting to see to what extent that can be done 
for $\epsilon $ in the range (\ref{0.18}). 

\par As in \cite{MeSj}, the \an{}ity assumptions seem to be quite essential 
at least in the case of fixed $\epsilon $. Indeed one is naturally led to 
work in modified Hilbert spaces defined by introducing microlocal 
exponential weights in the spirit of \cite{Sj6, HeSj, MeSj, HiSj}, and
there are closely related \fop{}s with complex phase some of which have
associated complex \ctf{}s that are $\epsilon $-perturbations of the
identity.  

\par The plan of the paper is the following:
\smallskip
\par\noindent In section \ref{SectionGe}, we make the geometrical work and 
construct invariant torii close to the real domain. This allows us to 
construct a complex \ctf{} which reduces $p$ to a function on the 
cotangent space on the standard 2-torus, which is \indep{} of the 
space-variables.
\smallskip

\par\noindent In section \ref{SectionSp} we perform further reductions for 
the whole \op{} and obtain a complete \asy{} description of all the \ev{}s of 
$P_\epsilon $ in a rectange of the form (\ref{0.17}). This is still  
somewhat formal, but 
\smallskip

\par\noindent in section \ref{SectionGr}, we introduce a global Grushin problem, 
and verify that the formal \ev{}s of the preceding section coincide 
modulo ${\cal O}(h^\infty )$ with the actual \ev{} in a rectangle 
(\ref{0.17}).


\section{Geometric reductions}\label{SectionGe}
\setcounter{equation}{0}

We use the notation and general set-up of the introduction. Thus let $p$ 
denote the semi-classical principal symbol of the unperturbed \op{}. As a warm-up we 
recall how the assumption (\ref{0.12}) follows from a seemingly weaker 
assumption. Thus replace (\ref{0.12}) by the assumption that 
for some $\alpha >0$, every point $\rho \in p^{-1}(]-\alpha
,\alpha [)$ belongs to a closed $H_p$-trajectory $\gamma (\rho )$ with 
period $T(\rho )>0$, depending
continously on $\rho$. Also assume $dp\ne 0$ on $\Gamma _0$. Then,
\smallskip
\par\noindent 1) If
$\gamma (\rho )\in \Gamma _E$ is the $T(\rho )$-periodic
$H_p$-trajectory passing through $\rho \in p^{-1}(E)$, then the
action $I(\gamma (\rho ))=\int_{\gamma (\rho )}\xi dx$ only
depends on $E$ but not on $\rho $.\smallskip
\par\noindent 2) We have the same conclusion for the period $T(\rho
)$ and hence (\ref{0.12}) holds.\smallskip

\par Indeed, consider first two trajectories $\gamma (\rho
_0),\gamma (\rho _1)$ and take an intermediate family $\gamma (\rho
_s)$, $0\le s\le 1$, depending continuously on $s$, so that the
union of the $\gamma (\rho _s)$ is a two-dimensional \mfld{}
$\Gamma \subset p^{-1}(E)$. Notice that ${\sigma _\vert }_{\Gamma }
=0$, since $H_p$ is tangent to $\Gamma $ and belongs to the radical
of the restriction of $\sigma $ to $p^{-1}(E)$. Hence by Stokes'
formula,
$$\int_{\gamma (\rho _1)}\xi dx-\int_{\gamma (\rho _0)}\xi dx=\int_\Gamma
 \sigma =0.$$
This shows 1). As for 2), let $\gamma _E\subset \Gamma _E$ be a
smooth family of $H_p$-periodic curves with period = $T(\gamma _E)$. Let $\Gamma =\cup_{E_0\le E\le E_1}\gamma
_E$ and let $\nu $ be a \vf{} on $\Gamma $, with $\nu (p)=1$. Let
$t$ be a multivalued time variable on $\Gamma $, so that
$H_pt=1$. Then we claim that
$${\sigma _\vert }_{\Gamma }=dp\wedge dt=d(pdt):$$
On the one hand, $\langle \sigma ,\nu \wedge H_p\rangle =\langle
dp,\nu \rangle =1$ and on the other hand
$$\langle dp\wedge dt,\nu \wedge H_p\rangle =\det\pmatrix{\langle
dp,\nu \rangle  &0\cr \langle dt,\nu \rangle  &\langle
dt,H_p\rangle }=1,$$
since the diagonal elements of the matrix are equal to 1, and the
claim follows.

\par By Stokes' formula,
\begin{eqnarray*}
\int_{\gamma (E_1)}\xi dx-\int_{\gamma (E_0)}\xi
dx=\int_\Gamma
\sigma =\int_\Gamma d(pdt)=
-\int _{\widetilde{\Gamma }}d(tdp)\\ =-\int _\alpha  t(\rho
)dp(\rho )+\int _\alpha  (t(\rho )+T(\rho ))dp=\int_\alpha T(\rho
)dp\\ =T(E_0)(E_1-E_0)+{\cal O}((E_1-E_0)^2),
\end{eqnarray*}
where $\widetilde{\Gamma }$ is the "rectangular domain" obtained by
placing a "cut" $\alpha $ from $\gamma (E_0)$ to $\gamma (E_1)$,
and we get the (well-known) formula,
$${d\over dE}I(\gamma (E))=T(\gamma (E)).$$
Since $I(\gamma (E))$ only depends on $E$ and not on the choice of
$\gamma (E)$, we get 2).

\par Let $p_\epsilon $ be as in the introduction, and let $q$ be defined 
in (\ref{0.15}). 
Let $G(x,\xi )$ be an analytic \fu{} defined in a \neigh{} of
$p^{-1}(0)$, such that
\ekv{Ge.5}
{
H_pG=q-\langle q\rangle ,
}
where we recall that $\langle q\rangle $ is the trajectory average, 
defined in (\ref{0.14}).

\par We will replace $T^*M$ by the new IR-\mfld{} $\Lambda
_{\epsilon G}=\exp (i\epsilon H_G)(T^*M)$ (defined in a
complex \neigh{} of $\Gamma _0$). Writing $\Lambda _{\epsilon G}\ni
(x,\xi )=\exp (i\epsilon H_G)(y,\eta )$, we use $\rho =(y,\eta )$ as
real symplectic coordinates on $\Lambda _{\epsilon G}$. By Taylor
expansion, we get
\eekv{Ge.6}
{p_\epsilon (\exp (iH_G(\rho )))=(p+i\epsilon q)(\exp (i\epsilon H_G(\rho
))+{\cal O}(\epsilon ^2)=}
{p(\rho )+i\epsilon (q-H_pG)(\rho )+{\cal O}(\epsilon ^2)=p+i\epsilon
\langle q\rangle +{\cal O}(\epsilon ^2). }

\par Recall the assumptions (\ref{0.16}), (\ref{0.21}), where we shall 
assume for simplicity that $F_0=0$. (This is no real restriction, since we 
can always replace $p_\epsilon $ by $p_\epsilon -i\epsilon F_0$.)
Since the Poisson bracket $\{ p,\Re\langle q\rangle \}$ is zero, we
see that every component of the set $\Lambda _{0,0}=\{ p=0,\, \langle
q\rangle
=0\}$ is a smooth Lagrangian torus. Assume for simplicity (as in the 
introduction), that we
only have one such component. Near this torus,
$p, \Re\langle q\rangle $ form an integrable system, so we can find a
real and analytic \canform{} $\kappa^{-1}$ from a \neigh{} of 
$\Lambda _{0,0}$ to a
\neigh{} of $\xi =0$ in $T^*{\bf T}^2$, so that
$p\circ \kappa $ and $\Re\langle q\rangle \circ \kappa $ (and hence also 
$\langle q\rangle \circ \kappa $ because of (\ref{0.16})) become functions
of $\xi $ only. Here ${\bf T}^2=({\bf R}/2\pi {\bf Z})^2$.

\par We can do this in the following way: Let $\Lambda _{E,F}$ be the
Lagrangian torus given by $p=E,\Re \langle q\rangle =F$, for $(E,F)\in
{\rm neigh\,}(0,{\bf R}^2)$. Let $\gamma _1(E,F)$ be the cycle in
$\Lambda_{E,F}$ corresponding to a closed $H_p$-trajectory with
minimal period, and let $\gamma _2(E,F)$ be a second cycle so that
$\gamma _1,\gamma _2$ form a fundamental system of cycles on the
torus
$\Lambda_{E,F}$. Necessarily $\gamma _2$ maps to the simple loop
given by
$\langle q\rangle =F$ in the abstract quotient manifold
$p^{-1}(E)/{\bf R}H_p$. Now it is classical (see Arnold [Ar]) that we
can find a real analytic canonical transformation $\kappa : {\rm
neigh\,}(\eta =0,T^*{\bf T}^2)\ni (y,\eta )\mapsto
(x,\xi )\in {\rm neigh\,}(\Lambda _{0,0},T^*{\bf R}^2)$, ${\bf
T}^2:=({\bf R}/2\pi {\bf Z})^2$ such that
\ekv{Ge.7}{\eta _j={1\over 2\pi }(\int _{\gamma _j(E,F)}\xi dx-\int_{\gamma
_j(0,0)}\xi dx),}
where $E,F$ depend on $(x,\xi )$ and are
determined by $(x,\xi )\in \Lambda _{E,F}$, i.e. by
$E=p(x,\xi ),F=\Re\langle q\rangle (x,\xi )$. (See also \cite{HiSj}.)

\par In the following we sometimes write $p$ instead of $p\circ
\kappa
$ and similarly for $\langle q\rangle $ (cf (\ref{0.16})):
$$p=p(\xi ),\,\, \langle q\rangle =\langle q\rangle (\xi ).$$
Then $H_p=\sum_1^2{\partial p\over \partial \xi _j}{\partial \over
\partial x_j}$. From (\ref{Ge.7}) and the discussion at the beginning of 
this section, we see that $p=p(\xi _1)$ in the new coordinates, so
\ekv{Ge.9}
{
H_p=c(\xi _1){\partial \over \partial x_1},\ p=p(\xi _1),\,
c={\partial p\over \partial \xi _1}\ne 0. }
The assumption (\ref{0.21}) implies:
\ekv{Ge.10}
{
{\partial p\over \partial \xi _1}\ne 0,\ {d\Re\langle q\rangle \over
\partial \xi _2}\ne 0. }
Thus
\ekv{Ge.11}
{p_\epsilon =p(\xi _1)+i\epsilon \langle q\rangle
(\xi )+r_\epsilon (x,\xi ),}
where $r_\epsilon ={\cal O}(\epsilon ^2)$ and $p$,
$\langle q\rangle $ satisfy (\ref{Ge.10}).

\par Now look for a "Lagrangian" torus $\Gamma $ in the
complex domain of the form
\ekv{Ge.12}
{
\xi =\phi '(x),\ x\in {\bf T}^2,
}
with $\phi $ grad-periodic (in the sense that $\nabla \phi $ is 
single-valued on ${\bf T}^2$) and complex-valued, $\phi '={\cal
O}(\widetilde{\epsilon })$, $\epsilon \ll \widetilde{\epsilon }\ll
1$, such that
${{p_\epsilon} _\vert }_\Gamma =0$. We get the eiconal equation
\ekv{Ge.13}
{p({\partial \phi \over \partial x_1})+i\epsilon \langle q \rangle
(\phi '_x)+r_\epsilon (x,\phi '_x)=0,}
where $r_\epsilon ={\cal O}(\epsilon ^2)$. Write $p(\xi _1)=
c\xi _1+{\cal O}(\xi _1^2)$,
$\langle q\rangle (\xi )=a\xi _1+b\xi _2+{\cal O}(\xi ^2)$, $c\in{\bf R}$, 
$c,\Re b\ne
0$ so that:
$$((c+ia\epsilon ){\partial \over \partial x_1}+i\epsilon
b{\partial \over \partial x_2})\phi +F_\epsilon (x,\phi '_x)=0,$$
where
$$F_\epsilon (x,\xi )={\cal O}(\epsilon ^2+\epsilon \xi ^2+\xi
_1^2).$$ For notational convenience, assume that $a=0$, $b,c=1$. Look
for
$\phi =\widetilde{\epsilon }\psi $, with $\psi '={\cal O}(1)$,
$\epsilon \ll\widetilde{\epsilon }\ll 1$. Then we get
\ekv{Ge.14}
{
({\partial \over \partial x_1}+i\epsilon {\partial \over \partial
x_2})\psi +G_{\epsilon ,\widetilde{\epsilon }}(x,\psi '_x)=0,}
with
$$G_{\epsilon ,\widetilde{\epsilon }}(x,\xi )={1\over
\widetilde{\epsilon }}F_\epsilon (x,\widetilde{\epsilon }\xi
)={\cal O}(\epsilon ({\epsilon \over \widetilde{\epsilon
}}+\widetilde{\epsilon })+\widetilde{\epsilon }\xi _1^2).$$

 \par Let $H^m({\bf T}^2)$ denote the standard Sobolev space of
order $m$. In the following estimates, $m>1$ is fixed.
 Using a standard result about non-linear functions of Sobolev
class functions (see \cite{AlGe}), we get
\smallskip
\par\noindent 1) If $(\epsilon ^{-1}\partial _{x_1},\partial
_{x_2})\psi ={\cal O}(1)$ in
$H^m$, then
$G_{\epsilon ,\tilde{\epsilon }}(x,\psi '_x)={\cal O}(\epsilon
({\epsilon \over \tilde{\epsilon }}+\widetilde{\epsilon }))$ in
$H^m$.\smallskip
\par\noindent 2) If $(\epsilon ^{-1}\partial _{x_1},\partial
_{x_2})\psi _j={\cal O}(1)$ in
$H^m$ for $j=0,1$, then,
$$\Vert [G_{\epsilon ,\widetilde{\epsilon }}(x,\psi _j')]_0^1\Vert
_{H^m}={\cal O}(\epsilon ({\epsilon \over
\widetilde{\epsilon }}+\widetilde{\epsilon }))\Vert ({1\over
\epsilon }\partial _{x_1},\partial _{x_2})(\psi _1-\psi _0)\Vert
_{H^m}.$$\smallskip

\par\noindent 3) If 
$$({\partial \over \partial x_1}+i\epsilon {\partial \over \partial 
x_2})u=v,\hbox{ with }u,v\hbox{ periodic,}$$
then 
$$\Vert (\epsilon ^{-1}\partial _{x_1},\partial _{x_2})u\Vert _{H^m}\le 
{C\over \epsilon }\Vert v\Vert _{H^m}.$$

\smallskip

\par We shall find solutions to (\ref{Ge.14}) that are grad-periodic \fu{}s of
the form
\ekv{Ge.15.5}
{\psi =\psi _{{\rm per}}+a(\epsilon x_1+ix_2)+b(\epsilon
x_1-ix_2),}
with a given complex constant $a={\cal O}(1)$, and where the
periodic \fu{} $\psi _{{\rm per}}$ and the complex constant $b$
will depend on $a$. If ${\cal F}u(k)$ denotes the \F{} \coef{} of
$u$ at $k$, we get the system:
\ekv{Ge.16}
{\cases{\displaystyle 2\epsilon b+{\cal F}(G_{\epsilon
,\widetilde{\epsilon }}(x,\psi _{{\rm per}}'+a(\epsilon
x_1+ix_2)'+b(\epsilon x_1-ix_2)'))(0)=0
\cr\displaystyle ({\partial \over \partial x_1}+i\epsilon
{\partial \over \partial x_2})\psi _{{\rm per}}+G_{\epsilon
,\widetilde{\epsilon }}(x,\psi _{{\rm per}}'+a(\epsilon
x_1+ix_2)'+b(\epsilon x_1-ix_2)')+2\epsilon b=0.}  }

\par We will find the solution as a limit of a sequence
$$\psi ^{(j)}=\psi _{{\rm per}}^{(j)}+a(\epsilon
x_1+ix_2)+b^{(j)}(\epsilon x_1-ix_2),\ j=0,1,2,...$$
with $\psi ^{(0)}=a(\epsilon x_1+ix_2)$, (and $b^{(0)}=0$, $\psi 
^{(0)}_{\rm per}=0$) where we impose
$$({\partial \over \partial x_1}+i\epsilon {\partial \over \partial
x_2})\psi ^{(j+1)}+G_{\epsilon ,\widetilde{\epsilon }}(x,{\psi
^{(j)}}')=0.$$
The last equation gives the following system analogous to (\ref{Ge.16}) 
that we label (${\rm S}_j)$:
$$
\cases{\displaystyle 2\epsilon b^{(j+1)}+{\cal
F}(G_{\epsilon ,\widetilde{\epsilon }}(x,{\psi _{{\rm
per}}^{(j)}}'+a(\epsilon x_1+ix_2)'+b^{(j)}(\epsilon
x_1-ix_2)'))(0) =0
\cr\displaystyle ({\partial \over \partial x_1}+i\epsilon
{\partial \over \partial x_2})\psi _{{\rm
per}}^{(j+1)}+G_{\epsilon ,\widetilde{\epsilon }}(x,{\psi^{(j)}
_{{\rm per}}}'+a(\epsilon x_1+ix_2)'+b^{(j)}(\epsilon
x_1-ix_2)')+2\epsilon b^{(j+1)}=0.}
$$

\par From ($\rm S_0$), and the facts 1), 3), we get
$$\vert b^{(1)}\vert ={\cal O}(1)({\epsilon \over
\widetilde{\epsilon }}+\widetilde{\epsilon }),$$
$$\Vert ({1\over \epsilon }\partial _{x_1},\partial
_{x_2})\psi _{{\rm per}}^{(1)}\Vert _{H^m}\le {\cal
O}(1)({\epsilon \over
\widetilde{\epsilon }}+\widetilde{\epsilon }).$$
For $j\ge 1$, we consider $({\rm S}_j)-({\rm S}_{j-1})$ and get, using also 2),
$$\vert b^{(j+1)}-b^{(j)}\vert \le {\cal O}(1)({\epsilon \over
\widetilde{\epsilon }}+\widetilde{\epsilon })(\Vert ({1\over
\epsilon }\partial _{x_1},\partial _{x_2}) (\psi_{{\rm per}}
^{(j)}-\psi_{{\rm per}} ^{(j-1)})\Vert _{H^m}+\vert
b^{(j)}-b^{(j-1)}\vert ),$$
\begin{eqnarray*}
&&\Vert ({1\over \epsilon }\partial _{x_1},\partial
_{x_2})(\psi_{{\rm per}} ^{(j+1)}-\psi_{{\rm per}} ^{(j)})\Vert
_{H^m}\le
\\
&&{\cal O}(1)({\epsilon
\over
\widetilde{\epsilon }}+\widetilde{\epsilon })(\Vert ({1\over
\epsilon }\partial _{x_1},\partial _{x_2})(\psi_{{\rm per}}
^{(j)}-\psi_{{\rm per}} ^{(j-1)})\Vert _{H^m}+\vert
b^{(j)}-b^{(j-1)}\vert ).
\end{eqnarray*}
This implies that
$$\vert b^{(j+1)}-b^{(j)}\vert +\Vert ({1\over \epsilon }\partial
_{x_1},\partial _{x_2}) (\psi_{{\rm per}} ^{(j+1)}-\psi_{{\rm per}}
^{(j)})\Vert _{H^m}\le ({\cal O}(1)({\epsilon \over
\widetilde{\epsilon }}+\widetilde{\epsilon }))^{j+1},$$ and since
$\epsilon \ll
\widetilde{\epsilon }\ll 1$, we see that the schema converges
towards a solution to (\ref{Ge.14}) of the form (\ref{Ge.15.5}), with
$a={\cal O}(1)$ (given), and
\ekv{Ge.17}{\vert b\vert +\Vert ({1\over \epsilon }\partial
_{x_1},\partial _{x_2})\psi _{{\rm per}}'\Vert _
{H^m}={\cal O}(1)({\epsilon \over
\widetilde{\epsilon }}+\widetilde{\epsilon }).}
For $\phi $ we then have
\ekv{Ge.18}
{
\phi =\phi _{{\rm per}}+\widetilde{\epsilon }a(\epsilon
x_1+ix_2)+\widetilde{\epsilon }b(\epsilon x_1-ix_2), }
with $\widetilde{\epsilon }a={\cal O}(\widetilde{\epsilon })$ given,
$\vert \widetilde{\epsilon }b\vert +\Vert (\epsilon ^{-1}\partial
_{x_1},\partial _{x_2})\phi _{{\rm per}}\Vert
_{H^m}={\cal O}(1)(\epsilon
+\widetilde{\epsilon }^2)$. In particular,
\ekv{Ge.19}
{
{\partial \phi \over \partial x_1}=a\widetilde{\epsilon }\epsilon
+{\cal O}(\epsilon (\epsilon +\widetilde{\epsilon }^2)),\ {\partial
\phi \over \partial x_2}=ia\widetilde{\epsilon }+{\cal O}(\epsilon
+\widetilde{\epsilon }^2). }

\par In this discussion $m$ is fixed and the estimates are
\uf{} \wrt{} $\epsilon $. Clearly $\phi $ only depends on the
choice of $\widetilde{\epsilon }a$ (with $m$ being fixed). As in
\cite{MeSj}, we see that $\phi $ depends \hol{}ally on
$\widetilde{\epsilon }a$, and extends \hol{}ally in $x$ to some
$(\epsilon,\widetilde{\epsilon }) $-dependent domain in such a way
that the dependence of $\widetilde{\epsilon }a$ is still \hol{}. In
the preceding constructions, everything works the same way, if we
replace ${\bf T}^2$ by ${\bf T}^2+iy$, $\vert y\vert <1/C$, so it
follows that $\phi $ extends in $x$ to a complex \neigh{} of the
real torus, which is \indep{} of $\epsilon ,\widetilde{\epsilon }$,
and that the preceding estimates remain valid here.

\par Write $\phi =\phi _a$, when $\widetilde{\epsilon }$ is
fixed. Let $\Gamma _\phi $: $\xi =\phi '(x),\,\, x\in {\bf T}^2$. Let
$I_j(\Gamma _\phi )$, $j=1,2$ be the corresponding actions \wrt{} $\xi
_1dx_1+\xi _2dx_2$. From
(\ref{Ge.17}), (\ref{Ge.18}) (or simply (\ref{Ge.19})) we get:
\eekv{Ge.20}
{
I_1(\Gamma _\phi )=2\pi \widetilde{\epsilon }\epsilon (a+b)=2\pi
\widetilde{\epsilon }\epsilon (a+{\cal O}({\epsilon \over
\widetilde{\epsilon }}+\widetilde{\epsilon })), }
{
I_2(\Gamma _\phi )=2\pi i\epsilon (a-b)=2\pi
\widetilde{\epsilon } (ia+{\cal O}({\epsilon \over
\widetilde{\epsilon }}+\widetilde{\epsilon })), }
We are interested in finding $a$ such that both actions are real.
This leads to
$${\rm Im\,}(a+{\cal O}({\epsilon \over \widetilde{\epsilon
}}+\widetilde{\epsilon }))=0,\ {\rm Im\,}(ia+{\cal O}({\epsilon
\over \widetilde{\epsilon }}+\widetilde{\epsilon }))=0,$$
i.e.
\ekv{Ge.21}
{
\cases{\displaystyle
{\rm Re\,}a+{\cal O}({\epsilon \over \widetilde{\epsilon
}}+\widetilde{\epsilon })=0\cr\displaystyle {\rm Im\,}a+{\cal
O}({\epsilon \over \widetilde{\epsilon }}+\widetilde{\epsilon
})=0.}
}
Here the ${\cal O}$-terms are real parts of \hol{} \fu{}s, so they
remain ${\cal O}(\widetilde{\epsilon }+\epsilon
/\widetilde{\epsilon })$ after derivation \wrt{} ${\rm Re\,}a$,
${\rm Im\,}a$. By the implicit \fu{} \th{}, we therefore have a
unique solution to (\ref{Ge.21}), which is
${\cal O}(\widetilde{\epsilon}+\epsilon /\widetilde{\epsilon })$,
and correspondingly $\widetilde{\epsilon }a={\cal O}(\epsilon
+\widetilde{\epsilon }^2)$. Recall that $\widetilde{\epsilon }a$
is independent of the choice of $\widetilde{\epsilon }$, so if we
take $\widetilde{\epsilon }=\sqrt{\epsilon }$, we get
\ekv{Ge.22}
{
\widetilde{\epsilon }a={\cal O}(\epsilon ).
}
For this particular $\phi $, we have
\ekv{Ge.23}
{
\partial _{x_1}\phi ={\cal O}(\epsilon ^2),\ \partial _{x_2}\phi
={\cal O}(\epsilon )\hbox{ in }H^m. }

\par If we do not make the simplifying assumption that ${\partial
p\over \partial \xi _1}(0)=1$, ${\partial \langle q\rangle \over
\partial \xi _1}(0)=0$,  ${\partial \langle q\rangle \over \partial
\xi _2}(0)=1$, then the earlier discussion tells us that we have
solutions of the type $\phi =\widetilde{\epsilon }\psi $, $\psi
=\psi _{{\rm per}}(x)+a\alpha (x)+b\beta (x)$, with  $$\alpha
(x)=\epsilon {\partial \langle q\rangle \over \partial \xi _2}(0)x_1
+ i {\partial (p+i\epsilon \langle q\rangle )\over \partial \xi
_1}(0)x_2,$$ $$\beta (x)=\epsilon {\partial \langle q\rangle \over
\partial \xi _2}(0)x_1 - i {\partial (p+i\epsilon \langle q\rangle
)\over \partial \xi _1}(0)x_2.$$ Observe that if we put,  $$Z:=
{\partial (p+i\epsilon \langle q\rangle )\over \partial \xi
_1}(0){\partial
\over \partial x_1}+{\partial i\epsilon \langle q\rangle \over
\partial \xi _2}(0){\partial \over \partial x_2},$$
then
$$Z\alpha =0,\ Z\beta
=2\epsilon {\partial (p+i\epsilon \langle q\rangle )\over \partial \xi
_1}(0){\partial \langle q\rangle \over \partial \xi _2}(0)\ne 0.$$

\par The earlier discussion goes through without any
changes. Especially, in the case $a=0$, the corresponding $\phi $ is
\indep{}  of $\widetilde{\epsilon }$.

\par Let now $\zeta $ vary in ${\rm neigh\,}(0,{\bf C}^2)$. Put $z(\zeta
)=p(\zeta _1)+i\epsilon \langle q\rangle (\zeta )$. Then the discussion
above can be applied with $p_\epsilon (x,\xi )$ replaced by
$$p_\epsilon (x,\zeta +\xi )-z(\zeta )=p(\zeta _1+\xi _1)-p(\zeta_1
)+i\epsilon (\langle q\rangle (\zeta +\xi )-\langle q\rangle (\zeta
))+{\cal O}(\epsilon ^2).$$ We get a solution to the eiconal equation
$$p_\epsilon (x,\zeta +\psi _x')-z(\zeta )=0$$ of the form $$\psi
(x,\zeta )=\psi _{\rm per}(x,\zeta )+b(\zeta )\beta (x,\zeta ),$$
where
$$\beta (x,\zeta )=\epsilon {\partial \langle q\rangle \over\partial
\xi _2}(\zeta )x_1-i{\partial (p+i\epsilon \langle q\rangle
)\over\partial \xi _1}(\zeta )x_2,$$ depending holomorphically on $\zeta
$. (So we choose $a=0$ in the earlier discussion, but compensate for
this by introducing a $\zeta $-dependence and even varying the
energy level $z(\zeta )$.)

\par As before, we get
\ekv{Ge.24}{ \Vert ({1\over \epsilon
}\partial _{x_1},\partial _{x_2})\psi _{\rm per} \Vert_{H^m}={\cal
O}(\epsilon ),\ |b|={\cal O}(\epsilon ), }  uniformly with respect to
$\zeta$. Moreover, since the problem depends holomorphically on
$\zeta $, it is easy to see (for instance by working in a space of
holomorphic functions of $\zeta $ with values in $H^m$)
that $\nabla _x\psi ,b$  depend holomorphically on $\zeta $. Notice
that
\ekv{Ge.25} {\widetilde{\psi }(x,\zeta ):=x\cdot \zeta
+\psi (x,\zeta )} solves the eiconal equation
\ekv{Ge.26}
{p_\epsilon (x,\partial _x\widetilde{\psi }(x,\zeta ))-z(\zeta )=0.}

\par Write $b(\zeta )\beta (x,\zeta )+x\cdot \zeta =x\cdot \eta $,
where $\eta (\zeta )$ depends \hol{}ally on $\zeta $ and satisfies
$$\eta _1(\zeta )=\zeta _1+{\cal O}(\epsilon ^2),\ \eta _2(\zeta
)=\zeta _2+{\cal O}(\epsilon ).$$
Let $\zeta (\eta )$ with $\zeta _1(\eta )=\eta _1+{\cal
O}(\epsilon ^2)$, $\zeta _2(\eta )=\eta _2+{\cal O}(\epsilon )$
denote the inverse. Then with $\phi _{\rm per}(x,\eta )=\psi _{\rm
per}(x,\zeta )$, we have
\ekv{Ge.27}
{\widetilde{\psi }(x,\zeta )=x\cdot \eta +\phi _{\rm per}(x,\eta
)=:\phi (x,\eta ),}
solving
\ekv{Ge.28}
{
p_\epsilon (x,\partial _x\phi (x,\eta ))-\widetilde{p}_\epsilon
(\eta )=0,\ \widetilde{p}_\epsilon (\eta )=z(\zeta (\eta ))=p(\eta
_1)+i\epsilon \langle q\rangle (\eta )+{\cal O}(\epsilon ^2), }
while (\ref{Ge.24}) gives
\ekv{Ge.29}
{
\vert ({1\over \epsilon }\partial _{x_1},\partial _{x_2})\phi _{\rm
per}(x,\eta )\vert ={\cal O}(\epsilon ), }
for $x$ in a fixed complex \neigh{} of ${\bf T}^2$ and as usual, we
get corresponding estimates for $\partial _x^\alpha \partial _\eta
^\beta \phi _{\rm per}(x,\eta )$ from the Cauchy inequalities. We
normalize the choice of $\phi _{\rm per}(x,\eta )$ by requiring
that $$\langle \phi _{\rm per}(\cdot ,\eta )\rangle ={1\over (2\pi
)^2}\int_{{\bf T}^2} \phi _{\rm per} (x,\eta )dx=0.$$
Then
\ekv{Ge.30}
{\kappa _\epsilon : (\phi _\eta '(x,\eta ),\eta )\mapsto (x,\phi
_x'(x,\eta ))}
maps a complex ($\epsilon $-\indep{}) \neigh{} of the zero section
of $T^*{\bf T}^2$ onto another \neigh{} of the same type (containing
an $\epsilon $-\indep{} \neigh{} of $\xi =0$). (\ref{Ge.28}) shows
that
\ekv{Ge.31}
{p_\epsilon \circ \kappa _\epsilon =\widetilde{p}_\epsilon .}
By construction, we also know that $\kappa _\epsilon $ conserves
actions along closed curves.

\par Using that $\phi '_\eta (x,\eta )=x+{\cal O}(\epsilon )$,
$\phi _x'(x,\eta )=\eta +{\cal O}(\epsilon ^2,\epsilon )$ together
with (\ref{Ge.29}), which also holds with $\phi _{\rm per}$ replaced
by its gradient, we see that
\ekv{Ge.32}
{
\kappa _\epsilon (y,\eta )=(y+{\cal O}(\epsilon );\eta_1 +{\cal
O}(\epsilon ^2),\eta _2+{\cal O}(\epsilon )). }
In particular, we have
\ekv{Ge.33}
{
\Im x={\cal O}(\epsilon ),\ \Im \xi_1 ={\cal O}(\epsilon ^2),\, \Im
\xi _2={\cal O}(\epsilon ), }
on the image of $T^*{\bf T}^2$. We can therefore represent $\kappa
_\epsilon (T^*{\bf T}^2)$ by
\ekv{Ge.34}
{
\Im x=G_\xi '(\Re (x,\xi )),\,\, \Im \xi =-G_x'(\Re (x,\xi )),
}
where $G$ is a smooth, a priori grad-periodic function which
satisfies,
\ekv{Ge.35}
{
\partial _\xi G,\partial _{x_2}G={\cal O}(\epsilon ),\ \partial
_{x_1}G={\cal O}(\epsilon ^2). }
Since $\kappa _\epsilon $ conserves actions, the actions along closed
cycles in $\kappa _\epsilon (T^*{\bf T}^2)$ are real and it follows
that $G$ is single-valued. We may assume that $G={\cal O}(\epsilon
)$. Let $\chi (\xi )$ be a standard cutoff around $\xi =0$ and let
$\widetilde{M}_\epsilon $ be given by
\eekv{Ge.36}
{
\Im x=\widetilde{G}'_\xi (\Re (x,\xi )),\, \Im \xi
=-\widetilde{G}'_x(\Re (x,\xi )),} 
{\hbox{where }\widetilde{G}(\Re (x,\xi
))=\chi (\Re \xi )G(\Re (x,\xi )). }
Then $\widetilde{M}_\epsilon $ is an IR-\mfld{} which coincides
with $T^*{\bf T}^2$ outside a (complex $\epsilon
$-indepen\-dent) \neigh{} of $\xi =0$. Moreover, we know that $\widetilde{M}
_\epsilon $ is an $\epsilon $-\pert{} of $T^*{\bf T}^2$, along which
we have
$$\Im \xi _1=-\chi (\Re \xi )G_{x_1}'(\Re (x,\xi ))={\cal
O}(\epsilon ^2). $$
It follows that outside the \neigh{} of $\xi =0$, where
$\widetilde{M}_\epsilon $ coincides with $\kappa _\epsilon (T^*{\bf
T}^2)$, we have
\ekv{Ge.37}
{
\vert \Re {{p_\epsilon }_\vert}_{\widetilde{M}_{\epsilon
}}\vert +{1\over \epsilon }\vert \Im {{p_\epsilon
}_\vert}_{\widetilde{M}_\epsilon }\vert \ge {1\over C}. }

\par Now recall the initial global situation, that we simplified
the original principal symbol by composing with $\exp i\epsilon
H_G$ for the \fu{} $G$ in (\ref{Ge.5}) and then further
by $\kappa $, introduced prior to (\ref{Ge.9}).

\par We introduce an IR-deformation $M_\epsilon $ of real phase space which is
an $\epsilon $-deformation, equal to real phase space away from
$\Gamma _0$, and equal to $\exp i\epsilon H_G\circ \kappa
(\widetilde{M}_\epsilon )$ near $\Gamma _0=p^{-1}(0)\cap T^*M$. Then we have achieved
the following:

\begin{prop}\label{PropGe.1} 
\par\noindent a) There exists an analytic real \canform{} $\kappa
_\epsilon :{\rm neigh\,}(\xi =0,T^*{\bf T}^2)\to {\rm
neigh\,}(\exp (i\epsilon H_G)(\Lambda _{0,0}),M_\epsilon )$, such that
\ekv{Ge.38}
{
p_\epsilon \circ \kappa _\epsilon =\widetilde{p}_\epsilon (\eta ),
}
where $\widetilde{p}_\epsilon $ is given in (\ref{Ge.28}).\smallskip
\par\noindent b) Away from the small \neigh{}, where (\ref{Ge.38})
holds, we have
\ekv{Ge.39}
{
\vert \Re {{p_\epsilon }_\vert}_{M_{\epsilon
}}\vert +{1\over \epsilon }\vert \Im {{p_\epsilon
}_\vert}_{M_\epsilon }\vert \ge {1\over C}. }
Here $p_\epsilon $ denotes the original principal
symbol of the perturbed \op{}.
\end{prop}

\par It is now clear that the main result of [MeSj] can be applied
to give the full \asy{}s for all \ev{}s of $P_\epsilon $ in a
domain $\vert \Re z\vert <1/{\cal O}(1)$, $\vert \Im z\vert
<\epsilon /{\cal O}(1)$, for $\epsilon >0$ small enough and for
$h<h(\epsilon )>0$ small enough depending on $\epsilon $. It is not
apriori clear however what kind of uniformity \wrt{} $\epsilon $ we
may have in this result. We shall employ quantum Birkhoff normal
forms in the next section and obtain a more uniform result, valid
for $\epsilon >h^\delta $ for any fixed $\delta \in ]0,{1\over 2}[$.

\section{Formal spectral \asy{}s}\label{SectionSp}
\setcounter{equation}{0}

As in \cite{HiSj} (see also \cite{MeSj}) we can implement $\kappa _\epsilon $ 
by an elliptic \fop{} $U=U_\epsilon :L^2_S({\bf T}^2)\to H(M_\epsilon )$ 
which is \ml{}ly defined from a \neigh{} of $\xi =0$ in $T^*{\bf T}^2$ to 
a \neigh{} of $\exp{i\epsilon H_G}(\Lambda _{0,0})$ in  
$M_\epsilon $. Here $S=(S_1, S_2)\in {\bf R}^2$, with $S_j=\int_{\gamma 
_j}\xi dx$, and $\gamma _j=\gamma _j(0,0)$
are introduced prior to (\ref{Ge.7}). $L_S({\bf T}^2)$ denotes the space of 
locally $L^2$-\fu{}s $u$ on ${\bf R}^2$ satisfying the Floquet periodicity 
condition:
\ekv{Sp.0}{u(x-\gamma )=e^{{i\gamma \over 2\pi }\cdot ({1\over h}S+{\pi 
\over 2}\alpha ^0)},\ \gamma \in (2\pi {\bf Z})^2,}
where $\alpha ^0=(\alpha _1^0,\alpha _2^0)\in{\bf Z}^2$ is a Maslov index. 
By abuse of notation, we still denote by $P_\epsilon $, the
conjugated \op{} $U_\epsilon ^{-1}P_\epsilon U_\epsilon $.

\par We have an analytic $h$-\pop{} $P_\epsilon $ on ${\bf T}^2$
(defined \ml{}ly near $\xi =0$), of order $0$ in $h$, with leading
symbol \indep{} of $x$:
\ekv{Sp.1}
{
p_\epsilon (\xi )=p(\xi _1)+i\epsilon \langle q\rangle (\xi )+{\cal
O}(\epsilon ^2), }
defined in a fixed complex \neigh{} of $\xi =0$ in $T^*{\bf T}^2$,
depending \hol{}ally on $\epsilon \in D(0,\epsilon _0)$. The full
symbol is
\ekv{Sp.2}
{
P_\epsilon (x,\xi ;h)=\sum_{j=0}^\infty  h^jp_j(x,\xi ,\epsilon ),
}
with $p_j(x,\xi ,\epsilon )$ \hol{} \wrt{} $(x,\xi )$ in a
$j$-\indep{} complex \neigh{} of $\xi =0$ and $C^\infty $ \wrt{} $\epsilon \in
[0,\epsilon _0[$, with $p_0(x,\xi ,\epsilon )=p_\epsilon (\xi )$.
Following the standard Birkhoff
normal form procedure, we shall remove the
$x$-dependence in the $p_j$ by means of conjugation by an elliptic
$h$-\pop{} of order $0$. Let $A$ be an $h$-\pop{} of order 0. Recall
that
$$e^APe^{-A}=e^{{\rm ad}_A}P=\sum {1\over k!}{\rm ad}_A^kP.$$
Let the full symbol of $A$ be of the form $\sum_{k=0}^\infty
h^ka_k$. Then on the \op{} level,
\begin{eqnarray*}e^APe^{-A}&=&\sum_{\ell =0}^\infty \sum_{k=0}^\infty \sum_{j_1=0}^\infty
...\sum_{j_k=0}^\infty  {1\over k!}h^{j_1+...+j_k+\ell +k}({1\over
h}{\rm ad}_{a_{j_1}})...({1\over h}{\rm ad}_{a_{j_k}})p_\ell\\
&=&\sum_{n=0}^\infty h^ns_n,\end{eqnarray*}
with $s_0=p_0$, $s_1={1\over i}H_{a_0}p_0+p_1=iH_{p_0}a_0+p_1$, ...,
$s_{n+1}=iH_{p_0}a_n+\widetilde{s}_{n+1}$,..., where
$\widetilde{s}_{n+1}$ only depends on $a_0,...,a_{n-1}$ and is the
sum of the \coef{}s for $h^{n+1}$ from the terms
$${1\over k!}h^{j_1+...+j_k+\ell +k}({1\over h}{\rm
ad}_{a_{j_1}})...({1\over h}{\rm ad}_{a_{j_k}})(p_\ell ),$$
with
$$j_1+...+j_k+\ell +k\le n+1,\ j_1,...,j_k<n, \hbox{ or } k=0,\,
\ell =n+1.$$

\par Notice that
$$H_{p_0}=H_{p_\epsilon }={\partial p(\xi _1)\over \partial \xi
_1}\partial _{x_1}+i\epsilon ({\partial \langle q\rangle \over
\partial \xi }+{\cal O} (\epsilon ))\cdot \partial _x,$$
and that we can solve
\ekv{Sp.3}
{H_{p_0}a=b(x,\xi )-\langle b(\cdot ,\xi )\rangle ,\ x\in {\bf T}^2,}
with $\Vert a\Vert _{H^{m+1}}\le {\cal O}(1)\epsilon ^{-1}\Vert
b\Vert _{H^m}$. As already noticed in the preceding section, the
same equation can be solved in a complex domain $\{ x\in{\bf T}^2;\,
\vert \Im x\vert <C_2\}$, and we get
\ekv{Sp.4}
{
\sup_{\vert {\mathrm Im\,} x\vert <C_2}\vert a(x,\xi )\vert \le
{C(C_1,C_2)\over \epsilon }\sup_{\vert {\mathrm Im\,} x\vert <C_1}\vert
b(x,\xi )\vert , }
if $C_1<C_2$. The shrinking of the domains in (\ref{Sp.4}) is not a
problem, since we can take a sequence of such domains with
$C_j\searrow C_\infty >0$.

\par By solving \e{}s of the type (\ref{Sp.3}), we can determine
$a_0,a_1,...$ successively, so that $s_j=s_j(\xi ,\epsilon )$ are
\indep{} of $x$.
Assume by induction that $\nabla a_j={\cal O}(\epsilon
^{-1-2j})$, for
$j\le n-1$ (in a complex domain, so that we have the same estimates
on the derivatives of $\nabla a_j$). Then the general term in
$\widetilde{s}_{n+1}$ is
$${\cal O}(1) \epsilon ^{-1-2j_1}...\epsilon
^{-1-2j_k}={\cal O}(1)({1\over \epsilon })^{2(j_1+..+j_k)+k}.$$ Here,
$$2(j_1+..+j_k)+k=2(j_1+..+j_k+k)-k\le 2(n+1-\ell )-k=2n+2-2\ell -k.$$ So
this quantity is ${\cal O}(1)({1\over \epsilon })^{2n}$ except possibly when
$2\ell +k<2$, i.e. when $k=\ell =0$ or when $k=1,\,
\ell =0$. In the first case we get the \coef{} for $h^{n+1}$ in
$p_0$ which is
$0$. In the second case, we get the \coef{} for $h^{n+1}$ in
$h^{j_1+1}({1\over h}{\rm ad}_{a_{j_1}})(p_0)$ with
$j_1<n$, which is ${\cal O}(1)({1\over \epsilon })^{1+2j_1}$. Here
$1+2j_1\le 2n$. Thus
$\widetilde{s}_{n+1}={\cal O}(\epsilon ^{-2n})$ (in a complex domain). We
can choose $a_n$ periodic, with
$iH_{p_0}a_n=-\widetilde{s}_{n+1}+\langle
\widetilde{s}_{n+1}(\cdot ,\xi )\rangle $ and with $\nabla
a_n={\cal O}(\epsilon ^{-1-2n})$. This completes the induction step and we
find $a_k$ with $\nabla a_k={\cal O}(\epsilon ^{-1-2k})$ in a fixed complex
\neigh{} of ${\bf T}^2\times \{ \xi=0\}$
such that if
$$A^{(N)}=\sum_{k=0}^{N-1}{h}^ka_k,$$
then
\ekv{Sp.5}{\widetilde{P}^{(N)}:=e^{A^{(N)}}P_\epsilon
e^{-A^{(N)}}=\sum_{n=0}^\infty
h^n\widetilde{p}_n^{(N)},}
where $\widetilde{p}_n^{(N)}(\xi ,\epsilon )={\cal O}(\epsilon ^{-2(n-1)_+})$ and
$\widetilde{p}_n^{(N)}=\widetilde{p}_n^{(\infty )}$ is \indep{} of
$x$ and $N$, for $n\le N$. From this we
get the following formal spectral result:
\begin{theo} \label{ThSp.1}  Under the assumptions above,
there exists a constant $C>0$ such that if $\delta >0$ is fixed and
$h^{{1\over 2}-\delta }<\epsilon <1/C$, and $0<h\le h(\delta )$
with $h(\delta )>0$ small enough, then in the
region
\ekv{Sp.6}
{\vert \Re z\vert <{1\over C},\ {\vert \Im z\vert \over \epsilon
}<{1\over C},}
$P$ has the following quasi-\ev{}s:
\ekv{Sp.7}
{
z_k\sim \sum_{n=0}^\infty  h^n\widetilde{p}_n^{(\infty )}(h(k-{S\over 2\pi 
h}-{\alpha ^0\over 4}),\epsilon ),\ k\in{\bf Z}^2.
}
Here, $S\in{\bf R}^2$, $\alpha ^0\in {\bf Z}^2$ were introduced in the 
beginning of this section, and $p_0^{(\infty )}=p_\epsilon $ is given in 
(\ref{Sp.1}).\end{theo}

\par We leave undefined, the notion of quasi-\ev{}, and interpret the above 
theorem as the formal consequence of the reductions above and the fact that 
the functions 
$$e_k(x)=e^{ix\cdot (k-{S\over 2\pi h}-{\alpha ^0\over 4})},\ k\in{\bf Z}^2,$$
form an orthonormal basis in $L^2_S({\bf T}^2)$.

\section{Justification via a global Grushin \pb{}.}\label{SectionGr}
\setcounter{equation}{0}

In this section we outline how \Th{} \ref{ThSp.1} actually gives all 
\ev{}s in the rectangle (\ref{Sp.6}). As in \cite{MeSj}, \cite{HiSj} we 
construct an auxiliary, so called Grushin \pb{}. Actually, this construction 
is identical with the one in \cite{HiSj}, so we shall only recall 
the main steps.

\par For $C>0$ \sufly{} large, let $I(C,\epsilon )$ (depending also on $h$) 
be the set of all $k\in{\bf Z}^2$, for which the values $z_k$ in 
(\ref{Sp.7}) belong to the rectangle (\ref{Sp.6}). Recall that $z_k$ 
correspond to the orthonormal \fy{} of \fu{}s $e_k$, defined after \Th{} 
\ref{ThSp.1}.

\par Let $\kappa _\epsilon ,M_\epsilon $ be as in Proposition, \ref{PropGe.1} 
and let $U_\epsilon $ be the \fop{} quantization of $\kappa _\epsilon $ 
introduced in the beginning of Section \ref{SectionSp}. With $A^{(N)}$ 
defined there, let $A$ be a natural \asy{} limit. Define
\ekv{Gr.1}
{R_+:H(M_\epsilon )\to {\bf C}^{I(C,\epsilon )},}
by 
\ekv{Gr.2}
{R_+u(k)=(e^AU_\epsilon ^{-1}u\vert e_k)_{L^2_S}.}
Notice that $R_+$ is a globally welldefined \op{} modulo some 
indetermination of norm ${\cal O}(h^\infty )$, since $e^AU_\epsilon ^{-1}u$ 
is \ml{}ly welldefined in a \neigh{} of the zero section in $T^*{\bf T}^2$. 
Similarly, we define $R_-:{\bf C}^{I(C,\epsilon )}\to H(M_\epsilon )$, by 
\ekv{Gr.3}
{R_-u_-=\sum_{k\in I(C,\epsilon )}u_-(k)U_\epsilon e^{-A}e_k.}

\par Then for $z$ in the rectangle (\ref{Sp.6}), with an increased value 
of $C$, the \pb{}
\ekv{Gr.4}{(P-z)u+R_-u_-=v,\ R_+u=v_+,}
has a unique solution $(u,u_-)\in H(M_\epsilon )\times {\bf 
C}^{I(C,\epsilon )}$ for every $(v,v_+)\in H(M_\epsilon )\times {\bf 
C}^{I(C,\epsilon )}$. (Here we assume for simplicity that $P$ is a \bdd{} 
\op{}, otherwise we would have to work with modifications of 
$H(M_\epsilon )$ of Sobolev type, depending on additional order \fu{}s. See 
the appendix in \cite{HiSj} for more details and further references.) We 
have the corresponding apriori estmate
\ekv{Gr.5}{\Vert u\Vert +\Vert u_-\Vert \le {C\over \epsilon }(\Vert v\Vert 
+\epsilon \Vert v_+\Vert ),}
and if we write the solution
\ekv{Gr.6}{\pmatrix{u\cr u_-}=\pmatrix{E &E_+\cr E_- &E_{-+}}\pmatrix{v\cr 
v_+},}
then modulo ${\cal O}(h^\infty )$, $E_{-+}$ is the diagonal matrix 
$((z-z_k)\delta _{j,k})$, where $z_k$ are given in (\ref{Sp.7}).

\par Recall from \cite{HiSj} that the verification of these facts consists 
of half-estimates away from $\Lambda _{0,0}$ and the exploitation near 
$\Lambda _{0,0}$ of the reduction to a translation invariant \op{} on ${\bf 
T}^2$ in 
the preceding section. Since the \ev{}s of $P$ in our rectangle are 
precisely the values $z$ for which $E_{-+}(z)$ is non invertible, we get
\begin{theo} \label{ThGr.1}  Under the assumptions of Theorem \ref{ThSp.1},
there exists a constant $C>0$ such that if $\delta >0$ is fixed and
$h^{{1\over 2}-\delta }<\epsilon <1/C$, and $0<h\le h(\delta )$
with $h(\delta )>0$ small enough, then in the
region
\ekv{Gr.7}
{\vert \Re z\vert <{1\over C},\ {\vert \Im z\vert \over \epsilon
}<{1\over C},}
the \ev{}s of $P$ are simple and given by 
\ekv{Gr.8}
{
z_k\sim \sum_{n=0}^\infty  h^n\widetilde{p}_n^{(\infty )}(h(k-{S\over 2\pi 
h}-{\alpha ^0\over 4}),\epsilon ),\ k\in{\bf Z}^2,
}
with one \ev{} for each $k$ such that $z_k$ belongs to (\ref{Gr.7}).
Here, $S\in{\bf R}^2$, $\alpha ^0\in {\bf Z}^2$ were introduced in the 
beginning of Section \ref{SectionSp}, and the $\widetilde{p}^{(\infty 
)}_n$ were constructed prior to Theorem \ref{ThSp.1}. Further, 
$p_0^{(\infty )}(\xi ,\epsilon )=p(\xi _1)+i\epsilon \langle q\rangle 
(\xi )+{\cal O}(\epsilon ^2)$.

\end{theo}

\section{Application to barrier top \res{}s.}\label{SectionBT}
\setcounter{equation}{0}

We extend the domain of validity of one of the results of section 7 in 
\cite{HiSj}, by using \Th{} \ref{ThGr.1} as the new ingredient. The 
discussion that follows will therefore only be a brief recollection of a 
part of Section 7 in \cite{HiSj}, and we refer to that work for more details. Let 
\ekv{BT.1}{P=-h^2\Delta +V(x),\ p(x,\xi )=\xi ^2+V(x),\ (x,\xi )\in T^*{\bf 
R}^2={\bf R}^4,}
satisfy the general conditions for defining \res{}s near the 
energy level $E_0>0$. Assume that $V(0)=0$, $\nabla V(0)=0$, $V''(0)<0$ and 
that $V$ is everywhere analytic. After a linear change of 
$x$-coordinates, we have near $x=0$:
\ekv{BT2}
{p(x,\xi )-E_0=\sum_1^2{\lambda 
_j\over 2}(\xi _j^2-x_j^2)+p_3(x)+p_4(x)+...,}
where $\lambda _j>0$ and $p_\nu $ is a homogeneous \pol{} of degree 
$\nu $. Also assume that (0,0) is the only trapped point for the $H_p$-flow 
on the real energy surface $p^{-1}(E_0)$.

We assume $\lambda =(\lambda _1,\lambda _2)$ fulfills the \res{} condition 
\ekv{BT.3}
{\lambda \cdot k=0,\hbox{ for some }0\ne k\in{\bf Z}^2.}
Somewhat roughly, the problem of determining the \res{}s near $E_0$ is then 
equivalent to determining the \ev{}s of $P-E_0$ near $0$, after the 
change of variables, $x=e^{i\pi /4}\widetilde{x}$ (and $\xi 
=e^{-i\pi /4}\widetilde{\xi }$) near $0$, and we get a new \op{} with symbol
\eekv{BT.4}{&&-i(p_2(\widetilde{x},\widetilde{\xi})+ie^{3\pi i/4}p_3(\widetilde{x})+
ie^{4\pi i/4}p_4(\widetilde{x})+...)=-iq(\widetilde{x},\widetilde{\xi}),} 
{&&p_2(\widetilde{x},\widetilde{\xi })=\sum_1^2{\lambda _j\over 
2}(\widetilde{\xi} 
_j^2+\widetilde{x}_j^2).}
Dropping the tildes for the new variables, we are then interested in \ev{}s 
$E$ of $Q=q(x,hD_x)$ with $\vert E\vert \sim \epsilon ^2$, $h^\delta <\epsilon \ll 
1$, $0<\delta <1/2$. Write $x=\epsilon y$, $\widetilde{h}=h/\epsilon ^2$. Then 
$hD_x=\epsilon \widetilde{h}D_y$ and 
$$\epsilon ^{-2}q(x,\xi )=\epsilon ^{-2}q(\epsilon (y,\eta ))=p_2(y,\eta 
)+i\epsilon e^{3\pi i/4}p_3(y)+{\cal O}(\epsilon ^2),$$
in a region $\vert (y,\eta )\vert ={\cal O}(1)$, where the corresponding \ef{}s are 
concentrated.

\par The \res{} condition (\ref{BT.3}) implies that the $H_{p_2}$-flow is 
periodic with a period $T>0$, \indep{} of the energy level. Using \Th{} 
\ref{ThGr.1} in the discussion of section 7 in \cite{HiSj}, we get the 
following variant of Proposition 7.1 of that paper:

\begin{prop}\label{PropBT.1}
Let 
$\langle{p_3}\rangle$ denote the average of $p_3$ along the trajectories 
of the Hamilton vector field of 
$p_2$ in {\rm (\ref{BT.4})}, and assume that $\langle{p_3}\rangle$ is not 
identically zero. Let $F_0 \in {\bf R}$ be a regular value of 
$\cos(3 \pi /4)\langle{p_3}\rangle$ restricted to $p_2^{-1}(1)$, and 
assume that $T$ is the minimal period of the 
$H_{p_2}$-trajectories in the torus $\Lambda_{1,F_0}$ given by 
$$
\Lambda_{1,F_0}: p_2 = 1, \cos \left(\frac{3\pi}{4}\right) 
\langle{p_3}\rangle = F_0. 
$$
Let $\epsilon $ satisfy 
\ekv{BT.5}
{h^\delta  \ll \eps \le \epsilon _0,\ 0<\epsilon _0\ll 1,\ 0<\delta 
<{1\over 4}.
}
Then for $z$ in the set 
$$
\left[1-\frac{1}{{\cal O}(1)}, 1+\frac{1}{{\cal O}(1)}\right] + i \eps 
\left[F_0 - \frac{1}{{\cal O}(1)}, 
F_0 + \frac{1}{{\cal O}(1)}\right], 
$$
the \res{}s of the form $E_0-i\epsilon ^2z$ are 
given by
$$
z = \widehat{P}\left(\widetilde{h}(k-\frac{\alpha}{4})-\frac{S}{2\pi}, 
\eps ; \widetilde{h}\right)+
{\cal O}(h^{\infty}),\quad \widetilde{h}=\frac{h}{\eps^2},\;\; k\in \z^2. 
$$
(with precisely one \res{} for every $k$).
Here $\widehat{P}\left(\xi,\epsilon ; 
\widetilde{h}\right)$ has an expansion as $\widetilde{h}\rightarrow 0$, 
$$
\widehat{P}\left(\xi,\eps; 
\widetilde{h}\right)\sim \sum_{n=0}^\infty 
\widetilde{h}^n\widetilde{p}_n^{(\infty )}(\xi ,\epsilon ),$$
where 
$$
\widetilde{p}_0(\xi ,\epsilon )=p_2(\xi )+i\epsilon e^{3\pi i/4}\langle 
p_3\rangle (\xi )+{\cal O}(\epsilon ^2),\ \widetilde{p}_j(\xi ,\epsilon 
)={\cal O}(\epsilon ^{-2(j-1)}),\,\, j\ge 1.
$$
The coordinates $\xi_1=\xi_1(E)$ and $\xi_2=\xi_2(E,F)$ are the 
norma\-lized ac\-tions of 
$$
\Lambda_{E,F}: p_2=E,\ \cos\left(\frac{3\pi}{4}\right) 
\langle{p_3}\rangle=F,
$$
for $E\in {\rm neigh}(1,\real)$, $F\in {\rm neigh}(F_0,\real)$, given by 
\begeq
\label{7.10}
\xi_j = \frac{1}{2\pi} \left(\int_{\gamma_j(E,F)} 
\eta\,dy-\int_{\gamma_j(1,F_0)} \eta \,dy\right), \quad j=1,2,
\endeq
with $\gamma_j(E,F)$ being fundamental cycles in $\Lambda_{E,F}$, such 
that $\gamma_1(E,F)$ corresponds to a closed $H_{p_2}$-trajectory of 
minimal 
period $T$. Furthermore, 
\begeq
\label{7.11}
S_j =\int_{\gamma_j(1,F_0)} \eta\,dy, \quad j=1,2,\,\,\, S=(S_1,S_2),  
\endeq
and $\alpha \in \z^2$ is fixed.  
\end{prop}

\par The interest of this result (as well as of \Th{} \ref{ThGr.1}) compared 
to the corresponding ones in 
\cite{HiSj} is that we can reach small but $h$-\indep{} values of 
$\epsilon $. On the other hand our method does not immediately seem to be 
able to hand as small values of $\epsilon $ as in \cite{HiSj}, and the 
results there give a desrciption of how the negative powers of $\epsilon $ 
appear in our estimates of the terms in the \asy{} expansion of the 
symbol $\widehat{P}(\xi ,\epsilon ;h)$.

\end{document}